\documentclass[11pt,reqno]{amsart}

\usepackage{amsmath}
\usepackage{amssymb}
\usepackage{amsthm}
\usepackage{amscd, amsfonts, mathrsfs, xcolor}
\usepackage{amsaddr}

\usepackage{cases}

\usepackage{verbatim}

\usepackage{epsf}
\usepackage[bookmarksnumbered,pdfpagelabels=true,plainpages=false,
colorlinks=true,
            linkcolor=black,citecolor=black,urlcolor=black]{hyperref}

\theoremstyle{plain}
\newtheorem{theorem}{Theorem}[section]

\theoremstyle{definition}
\newtheorem{remark}[theorem]{Remark}

\numberwithin{equation}{section}
\allowdisplaybreaks

\newcommand{\al}{\alpha}
\newcommand{\be}{\beta}

\newcommand{\R}{\mathbb R}

\newcommand{\vol}{\operatorname{Vol}}

\newcommand{\dive}{\operatorname{div}}
\newcommand{\grad}{\operatorname{grad}}
\newcommand{\curl}{\operatorname{curl}}
\newcommand{\Ric}{\operatorname{Ric}}
\newcommand{\Def}{\operatorname{Def}}

\def\beq{\begin{equation}}
\def\ee{\end{equation}}
\newcommand{\NS}[1]{\N_{#1}}

\newcommand{\dd} {\mathrm{d}}

\usepackage{centernot}

\title[Navier-Stokes on manifolds]{The formulation of the Navier-Stokes equations on Riemannian manifolds}

\author[Chan]{\vspace{0.0cm} Chi Hin Chan}
\address{\vspace{-0.25cm} Department of Applied Mathematics\\
National Chiao Tung University\\ 1001 Ta Hsueh Road, Hsinchu, Taiwan 30010, ROC}
\email{cchan@math.nctu.edu.tw}

\author[Czubak]{\vspace{0.0cm} Magdalena Czubak}
\address{\vspace{-0.25cm} Department of Mathematics\\
University of Colorado Boulder\\ Campus Box 395, Boulder, CO, 80309, USA}
\email{magda.czubak@colorado.edu}

\author[Disconzi]{\vspace{-0.0cm} Marcelo M. Disconzi}
\address{\vspace{-0.25cm} Department of Mathematics\\
Vanderbilt University \\
1326 Stevenson Center,  Nashville, TN, 37240, USA}
\email{marcelo.disconzi@vanderbilt.edu}

\DeclareMathOperator{\N}{N-S}

\begin{document}

\begin{abstract}
We consider the generalization of the Navier-Stokes equation from $\R^n$ to the Riemannian manifolds.
There are inequivalent  formulations of the Navier-Stokes
equation on manifolds due to the different possibilities for the Laplacian operator acting on vector fields
on a Riemannian manifold.  We present several distinct arguments that
indicate that the form of the equations proposed by Ebin and Marsden in 1970 should be
adopted as the correct generalization of the Navier-Stokes to the Riemannian manifolds.
\end{abstract}

\maketitle

\tableofcontents

\section{Introduction \label{intro}}
The Navier Stokes equations are one of the fundamental equations of fluid
mechanics. They play an important role in aerodynamics, geophysics, meteorology,
and
engineering.  On $\R^{n}$ the equations are given by
\beq\tag{$\NS{\R^{n}}$}
\begin{split}
\partial_{t}v -\Delta v + \mu v\cdot \nabla v + \nabla p & = 0 , \\
\dive v& = 0 ,
\end{split}
\label{NSR}
\ee
where $v: \R^{n+1}\rightarrow \R^n$ is the velocity of the fluid, $p:
\R^{n+1}\rightarrow \R$ is the pressure, $\dive v = 0$ means the fluid is
incompressible, and $\mu$ is the viscosity.

Following the seminal work of Arnold \cite{ArnoldFluidsGeometric},
the study of the inviscid version of equations \eqref{NSR}, namely the Euler equations, in the setting of Riemannian manifolds,
has spurred a great deal of activity
and interplay between analysis and geometry.
The reader is referred to
the monograph \cite{ArnoldTopological} and references therein for an overview of the subject.
 
To the best of our knowledge, the first paper to present a systematic analysis of the Navier-Stokes
equations on the Riemannian manifolds is the work of Ebin and Marsden \cite{EbinMarsden}.  That work has been followed by a number of works in
\cite{AvezBamberger, CP_92, CC10, CCL15, CzubakChanRemarks, MR3119662, DindosMitrea,  Illin, IFilatov, MR3005957,  Lichtenfelz,   Mazzucato2003, MitreaTaylor, Nagasawa97, Nagasawa99,  Pierfelice, Priebe, Taylor_Morrey, TaylorPDE3, TemamWang, QSZhang}.

When one moves from the Euclidean setting to Riemannian manifolds, the first
question is
how to write the equations.  As numerous as the above works are, they do not all employ the same set of equations.   Hence it seems natural to inquire what the correct form of the equations should be.

 For the case of other important equations in physics and engineering, such as Euler or Maxwell equations,
the passage from $\R^n$ to a Riemannian manifold is more or less straightforward. This is because
such equations are obtained as critical points of an action functional that can be naturally
defined on a manifold, such as, for instance, the total energy in the case of the Euler equations.
The Navier-Stokes equations, however, do not come from an action, and thus it is not immediately
clear how to define them on a manifold.

One could argue that the equation
should be generalized directly from \eqref{NSR} upon interpreting each term by its corresponding
analogue on manifolds. For instance, $v \cdot \nabla v$ is simply the directional
derivative of $v$ in the direction of $v$, thus it should be generalized to
$\nabla_v v$, where $\nabla$ is the Levi-Civita connection associated with the metric.  Then we face the question how to interpret the viscosity operator, namely the Laplacian, $-\Delta$.
%

The first object that might come to mind upon hearing the words
Riemannian manifold and Laplacian is the Laplace-Beltrami operator
\[
 \dive \grad=\frac{1}{\sqrt{|g|}}\partial_{x^i}(\sqrt{|g|}g^{ij}\partial_{x^j}),
\]
where $|g|$ denotes the determinant of the metric $g=(g_{ij})$ written in local
coordinates $\{x^i\}$. However,
this operator acts on scalar valued functions. Meanwhile, a solution of the
Navier-Stokes equation is a vector field.  Hence, the Laplace-Beltrami operator cannot be applied here. Moreover, since we are considering vector fields, there are actually several
candidates for the choice of the Laplacian on a Riemannian manifold.

To begin with, if one likes to think about
the Laplacian as the $\dive \grad$ operator, then the natural generalization of
the Laplace-Beltrami operator to vector fields is the Bochner Laplacian
\[
 \dive \nabla=-\nabla^\ast\nabla,
\]
where $\nabla$ denotes the Levi-Civita connection on $(M,g)$, and $\nabla^\ast$
is the adjoint operator associated to $\nabla$. In local coordinates, we can
write
\[
 \dive \nabla= \nabla^i\nabla_i=g^{ij}\nabla_j \nabla_i.
\]
When $v$ is a vector field on $\R^3$, the Laplacian of $v$ can be expressed as
\[
 (\grad \dive  - \curl\curl )v.
\]
Now, when $v$ is a vector field on a Riemannian manifold, then using the metric,
we can obtain a unique $1-$form $\alpha$ associated to $v$. (In local coordinates, if $v=v^i \partial_{x^i}$, then $\alpha=g_{ij}v^j dx^i$.)  The analog of the
above expression becomes
\[
 -(dd^\ast + d^\ast d)\alpha,
\]
which is the Hodge Laplacian acting on $\alpha$. Here $d^*$ is the formal adjoint of the exterior
differential operator $d$.

The Bochner Laplacian and the
Hodge Laplacian do coincide on $\R^n$, but in general, they are not the same.
They are related by the Bochner-Weitzenb\"ock formula (see, e.g., \cite[p. 561]{TaylorPDE3})
\begin{gather}
\nabla^\ast \nabla=dd^\ast +d^\ast d -\Ric,
\label{BW_formula}
\end{gather}
where $\Ric$ is the Ricci curvature.  So in the particular case of $\R^n$, where $\Ric\equiv 0$, we can see that these
operators do indeed coincide.

There is another operator that can be considered on a Riemannian manifold.  In
fact, in their 1970 article \cite{EbinMarsden}, Ebin and Marsden indicated that
when
writing the Navier-Stokes equation on an Einstein
manifold, the ordinary Laplacian should be replaced by the following operator
\[
 2\Def^\ast \Def,
\]
where $\Def$ is the deformation tensor, and $\Def^\ast$ is its adjoint (see section \ref{section_deformation}).  The
deformation tensor can be thought of as a symmetrization of the connection, and
in the coordinates we can write it as
\begin{gather}
 {(\Def v)}_{ij}=\frac 12 (\nabla_i v_j +\nabla_j v_i).
\label{deformation}
\end{gather}
Since  $2\Def^\ast \Def=2\dive \Def$, a direct computation using  (\ref{BW_formula}) and a Ricci identity gives (the details of this computation can be found in \cite[p. 562]{TaylorPDE3})
\beq\label{divdef}
 2\Def^\ast \Def=dd^\ast +d^\ast d+dd^\ast -2\Ric.
\ee
Then adopting \eqref{divdef} and using that $v$ is divergence free, so $d^\ast v=0$,
the Navier-Stokes equations read
\beq\tag{$\NS{\text{Riem}}$}
\begin{split}
\partial_{t}v  +  \nabla_v v  -\mu\Delta_H v - 2\mu \Ric(v) +  d p & = 0 , \\
d^\ast v& = 0 ,
\end{split}
\label{NS}
\ee
where $\Delta_H = - (dd^\ast + d^\ast d)$ and $(\Ric(v))_i = \Ric_{ij} v^j$.  (Here, we identify the vector field $v=v^i \partial_{x^i}$, with the $1$-form $\alpha=g_{ij}v^j dx^i$, and denote both simply by $v$.)

We note that while the article \cite{EbinMarsden} states that the ``correct" viscosity operator on a manifold should come from the deformation tensor, the article itself uses the Hodge Laplacian, and the deformation tensor is addressed only at the end, in the ``Note Added in Proof."  The follow up article \cite{CP_92} also uses the Hodge Laplacian, and so do \cite{AvezBamberger, Illin, IFilatov, Lichtenfelz, TemamWang}.  In these articles, the domain considered is that of a compact manifold, and often of a sphere, so heuristically speaking the lack of the Ricci term might not make a difference.  The exception to that is \cite{Lichtenfelz}, which considers a non-compact manifold, and in particular illustrates why it is only a heuristic whether the $\Ric$ term might not make a difference.  Indeed, in \cite{Lichtenfelz} the author shows the Navier-Stokes equation in 3D with a Hodge Laplacian will exhibit non-unique Leray-Hopf solutions.  It is not clear that this proof can go through on a 3D manifold if the deformation tensor is used.

Starting with the work of Taylor \cite{Taylor_Morrey}, the following authors have followed Ebin and Marsden in using the deformation tensor, even if for several of these works they worked either on a compact manifold or a compact subdomain: \cite{CC10, CCL15, CzubakChanRemarks, MR3119662, DindosMitrea, MR3005957,  Mazzucato2003, MitreaTaylor,  Nagasawa97, Nagasawa99, Pierfelice, Priebe, TaylorPDE3,  QSZhang}.

The purpose of this article is to provide further evidence of why one should use
the deformation tensor when studying the Navier-Stokes equations on Riemannian
manifolds, i.e., that equation (\ref{NS}) should be adopted.
We also would like to extend the discussion from Einstein manifolds to all general Riemannian manifolds.
 
We provide the following distinct arguments. The first is based on an energy estimate.
We show that if the Hodge-Laplacian is adopted, then a priori energy estimate for the solutions to the Stokes equation
is not possible.  Moreover, we give a counterexample to the energy equality in this setting. We begin with the Stokes equation, because as it is the linear version of the Navier-Stokes equation, it can be viewed as the starting point for the study of the PDE theory for the system.  We then adapt the counterexample to the Navier-Stokes equations.  Since energy estimates are the cornerstone of the existence theory
for the solutions for \eqref{NSR}, it is important to seek generalizations to manifolds that preserve
such tool, and at the same time lead to a consistent theory starting from the linear to the nonlinear equation. This is done in section \ref{section_energy}.

In section \ref{section_restriction}, 
we show that the desired viscosity operator can be obtained 
by considering the restriction of (\ref{NSR}) in 3D
to a sphere.

The third argument is based on the non-relativistic limit of the relativistic
Navier-Stokes equations. While the relativistic formulation of the Navier-Stokes equations
is also open to debate, the known proposals in the literature all lead
to same equations based on (\ref{deformation}) in the non-relativistic limit. This is discussed in
section \ref{section_nr_limit}.

We begin by taking a closer look at the deformation tensor, section \ref{section_deformation}, where we also analyze more conceptual arguments that come up in the derivation of the Navier-Stokes equations that go back to the aforementioned
work of Ebin and Marsden and also Serrin \cite{Serrin}.

\begin{remark}
It is interesting to notice that in a stochastic setting, the issue of the correct formulation of viscous
fluids on manifolds does not seem to arise. This is because in the stochastic case
it is possible to write appropriate functionals whose critical points give the equations of motion
\cite{MR2787078,EyinkActionPrinciple}.
\end{remark}

\section{The deformation tensor\label{section_deformation}}
We introduced the deformation tensor as
\begin{gather*}
 {(\Def v)}_{ij}=\frac 12 (\nabla_i v_j +\nabla_j v_i),
\end{gather*}
where $\nabla_i$ denoted the covariant derivatives.  On the Euclidean space they reduce to just regular derivatives.  Moreover the deformation tensor can make an explicit appearance in the Euclidean Navier Stokes equation if we recall the equation's meaning.

 To see that, the Navier Stokes equation \label{NS-R} is the equation of conservation of momentum for an incompressible fluid.  As mathematicians, we are used to seeing it as written in \eqref{NSR}, but engineers often write it as
 \beq\label{NSe}
 \partial_{t}v+ v\cdot \nabla v=  -\nabla p + \mu\Delta v,
 \ee
 or in a more general way
 \beq\label{NSe2}
 \partial_{t}v+ v\cdot \nabla v=f + \dive T,
 \ee
which we now explain.  The equation \eqref{NSe2} might be viewed more natural physically, because it explicitly breaks down the equation in the parts written in the conservation of momentum, equivalently in Newton's 2nd Law, well-known as
\[
F=ma,
\]
i.e., force equals mass times acceleration.  The left hand side in \eqref{NSe2} comes from $ma$, and the right hand side denotes all the forces acting on the fluid.  These consist of volume and surface forces.  Volume forces act on all elements of the volume of a continuum, and gravity is an example of a volume force.  We can denote the volume force per unit mass of the fluid by vector valued function $f$.

Next, the surface forces are what eventually can produce the deformation tensor, which in turn gives us the Laplacian.  First, the surface force acts on a surface element to which we assume we can assign a unit normal $n$.  Then, the surface force, as a force is also a vector, and its i'th component can be written as $T_{ij}n_j$, where $T$ is the stress tensor (see for example \cite{Batchelor}, and where we sum over repeated indices).  So if we consider a part of a fluid with volume $V$ and enclosed by a surface $S$, the total surface force acting on $S$ is given by
\[
\int_S T_{ij}n_j dS= \int_V \partial_j T_{ij} dV
\]
where the equality holds by the divergence theorem. So this is how we get \eqref{NSe2} for fluids with constant density, but what constitutes the stress tensor $T$?  For fluid at rest or also for perfect fluid, only normal stresses are exerted, and we have $T_{ij}=-p \delta_{ij}$, where $p$ is the pressure, $\delta_{ij}$ is the kronecker delta, and gives the familiar $\nabla p$ in the equation.  Now, for fluids in motion or non-perfect fluids, we also have tangential stresses, which for isotropic fluids arrive with the additional term in $T_{ij}$, which can be shown to be
\[
2\mu D_{ij}-\frac 23 \mu \dive v \delta_{ij},
\]
where $D_{ij}=(\Def v)_{ij}$ is the deformation tensor, $\frac 12(\partial_jv_i+\partial_i v_j)$, and $\mu$ is the viscosity coefficient  (again see \cite{Batchelor}). Then for an incompressible fluid we are just left with
\[
\mu (\partial_jv_i+\partial_i v_j),
\]
and after we take divergence as in \eqref{NSe2} we get exactly
\[
\mu \Delta v.
\]

Now, the first in-depth study of the Navier-Stokes
equations on the Riemannian manifolds was carried out by Ebin and Marsden \cite{EbinMarsden}.
There, the authors point out that the deformation tensor should be adopted when
writing the equations on manifolds. They further assume that the manifold was Einstein,
stressing that the physical assumptions in the derivation of the Navier-Stokes equations
may not be satisfied when the manifold is not Einstein.

 Ebin and Marsden refer to
the derivation of the Navier-Stokes equations given by Serrin \cite{Serrin}.
There, following Stokes, it is assumed that the stress-tensor $T$ satisfies the following properties (see section 59 of \cite{Serrin}):

1.  $T$ is a continuous function of the deformation tensor $ \Def$,
and is independent of all other kinematic quantities.

2. $T$ does not depend explicitly on the spatial position (spatial homogeneity).

3. There is no preferred direction in space (isotropy).

4. When $\Def = 0$, $T$ reduces to $-pI$, where $p$ is the pressure and $I$ the identity matrix.

It is not clear how one would assure 2 and 3 on a general manifold.
In fact, in his discussion of the Navier-Stokes in curvilinear coordinates (section 13), Serrin
remarks that, on a general manifold, it is not ``evident
how to formulate the principle of conservation of momentum."
However, right after such comments, he argues that
``there seems to be no valid objection to taking Eq. (12.3) as a postulate."
Equation (12.3) in question involves the divergence of the stress tensor,
and it gives (\ref{NS}) if $T$ satisfies 1 and 4 above along with the remaining
assumption used by Serrin in the derivation of the equations (such as, for instance,
that $T$ is linear in $\Def$, see again section 59 of \cite{Serrin}).

The above considerations also highlight a point that is sometimes obscured when the equations
are written as in \eqref{NSR}, namely, that it is the stress-tensor and the deformation
tensor, and not the Laplacian, that are of direct physical significance in the modeling that leads to
\eqref{NSR}. This naturally suggests that $T$ and $\Def$ should again be the primary objects
one considers on a manifold.

We finish this section remarking that there is yet further evidence that (\ref{NS}) should be adopted,
which is when boundary conditions are introduced. Indeed, Shkoller has showed in
\cite{ShkollerAveraged} that it is the formulation of the equations in terms of $\Def$ that is naturally
associated with the Dirchlet and Neumann boundary conditions for the fluid.

\section{Counterexample to an energy estimate\label{section_energy}}

We motivate the discussion by looking at the Euclidean setting. Let $\Omega$ to be
a bounded domain with smooth boundary in $\mathbb{R}^3$, and consider the
Cauchy problem for the following linear Stokes equation on $[0,T]\times \Omega$,
\begin{equation}\label{linearStokes}
\left\{
\begin{split}
\partial_t v - \mu \Delta v + \nabla P  & = f, \\
 \dive v & = 0 ,
\end{split}
 \right.
\end{equation}
where the constant $\mu > 0$ denotes the viscosity coefficient.
One usually employs the following function space in the classical theory of
weak solutions to
(\ref{linearStokes}):
\begin{equation}
V(\Omega) = \overline{\Lambda^1_{c , \sigma}(\Omega) }^{\|\cdot
\|_{H^1(\Omega)}} ,
\nonumber
\end{equation}
where $\Lambda_{c , \sigma }^1(\Omega)$ stands for the space of all smooth,
divergence free, compactly supported vector fields on $\Omega$. Since $\Omega$ is bounded, we use the
typical convention $\|\phi\|_{H^1(\Omega)} = \|\nabla \phi\|_{L^2(\Omega)}$ for
the $H^1$-norm. \\
To specify the admissible class of finite energy, divergence free initial data
for the Cauchy problem associated to \eqref{linearStokes}, we use
\begin{equation}
\textbf{H}(\Omega) =
\overline{\Lambda_{c,\sigma}^1(\Omega)}^{\|\cdot\|_{L^2(\Omega)}}.
\nonumber
\end{equation}
With an initial datum $v_0 \in \textbf{H}(\Omega)$ and an external force $f \in
L^2(0,T ;(V(\Omega))')$, a weak solution to \eqref{linearStokes} in $[0,T]\times
\Omega$ which arises from $v_0$ is an element
\begin{equation}
v \in C^0([0,T] ; \textbf{H}(\Omega)) \cap L^2(0,T ; V (\Omega)),
\nonumber
\end{equation}
satisfying the following properties:
\begin{itemize}
\item [A)] $\partial_t v \in L^2(0,T ; (V(\Omega))')$.
\item [B)]
\begin{equation}\label{WeaklinearStokes}
\big < \partial_t v (t) , \phi  \big >_{(V(\Omega))'\otimes V(\Omega)} +
\mu \int_{\Omega} \nabla v(t) : \nabla \phi =
\big < f(t) , \phi  \big >_{(V(\Omega))'\otimes V(\Omega)}
\end{equation}
holds for all $\phi \in V(\Omega )$, and for almost every $t\in [0,T]$.
\item [C)] $v(0) = v_0$.
\end{itemize}
Standard theory in the Euclidean setting (see for instance
\cite[Ch. 2]{SereginBook}) ensures the existence and uniqueness of such a weak solution
$v$ to the Cauchy problem of \eqref{linearStokes} with any prescribed initial
datum $v_0 \in H(\Omega )$ and external force $f \in L^2(0,T; (V(\Omega))' )$. We
would like to point out that the heart of the matter of such an existence and
uniqueness theory is an a priori estimate for the
quantity
\begin{equation}
\|v\|_{L^{\infty} (0,T ; \textbf{H}(\Omega ))}^2 + \|v \|_{L^2(0,T ; V(\Omega
))}^2
\nonumber
\end{equation}
in terms of $\|v_0\|_{L^2(\Omega)}$ and $\|f\|_{L^2(0,T ; (V(\Omega))')}$.
The derivation of such apriori estimate basically proceeds as follows. Take
$\phi = v (t)$ in \eqref{WeaklinearStokes} to obtain
\begin{equation}\label{Younginequality}
\begin{split}
\frac{1}{2} \frac{\dd}{\dd t} \|v (t)\|_{L^2(\Omega)}^2 + \mu \|\nabla v
(t)\|_{L^2(\Omega)}^2 &= \big < f(t) , v(t) \big >_{(V(\Omega))'\otimes V(\Omega)}
\\
& \leq \|f(t)\|_{(V(\Omega))'} \, \| v(t)\|_{V(\Omega)} \\
& = \|f(t)\|_{(V(\Omega))'} \, \|\nabla v(t)\|_{L^2(\Omega)} \\
& \leq \frac{1}{2\mu} \|f(t)\|_{(V(\Omega))'}^2 + \frac{\mu}{2} \|\nabla
v(t)\|_{L^2(\Omega)}^2,
\end{split}
\end{equation}
from which it follows that
\begin{equation}\label{Gronwall}
\frac{\dd}{\dd t} \|v (t)\|_{L^2(\Omega)}^2 + \mu \|\nabla v (t)\|_{L^2(\Omega)}^2
\leq \frac{1}{\mu} \|f(t)\|_{(V(\Omega))'}^2.
\end{equation}
So, integrating  \eqref{Gronwall} in $t$
immediately leads to the following
a priori estimate,
\begin{equation}\label{Euclideanapriori}
\|v\|_{L^{\infty}(0,T ; L^2(\Omega ))}^2 +  \mu \int_0^T
\int_{\R^2} \big | \nabla v (\tau )\big |^2
\dd x \dd \tau
 \leq
\|v_0\|_{L^2(\Omega)}^2 + \frac{1}{\mu} \|f\|_{L^2(0,T; (V(\Omega))')}^2.
\end{equation}
The derivation of \eqref{Euclideanapriori} by means of the Cauchy's inequality
type estimate as being done in \eqref{Younginequality} is indeed the backbone
that supports the classical existence and uniqueness theory for weak solutions
to the Cauchy problem \eqref{linearStokes}. However, the key point of the
estimate \eqref{Younginequality}   is the presence of the term $\mu \|\nabla
v(t)\|_{L^2(\Omega)}^2$ on the left-hand side of \eqref{Younginequality}, which
enables one to absorb the extra term $\frac{\mu}{2} \|\nabla v(t)\|_{L^2(\Omega)}^2$
on the right-hand side of \eqref{Younginequality}. 

These remarks give us some hints about what could possibly go
wrong in writing the analogous linear Stokes system in the general setting of a
Riemannian manifold, with the Hodge Laplacian $(-\Delta_H)$ as the choice of
the elliptic operator representing viscosity effect within the Navier-Stokes
flows.

To see that, suppose that we insist on working with the following version of
the linear Stokes equation on the hyperbolic space $\mathbb{H}(-a^2)$ with
constant sectional curvature $-a^2$.
\begin{equation}\label{fakeequation}
\left\{
\begin{split}
\partial_t v +\mu (-\Delta_H) v + d P & = f , \\
\dd^* v & = 0.
\end{split}
\right.
\end{equation}
We take as a norm $\|v\|_{H^1(\mathbb{H}^2(-a^2))} =
\|\nabla v\|_{L^2(\mathbb{H}^2(-a^2))}$, where $\nabla$ is the Levi-Civita
connection on $\mathbb{H}^2(-a^2)$ (see \cite[Estimate (2.2)]{CCL15}, which guarantees this is a norm). Then, it is natural to look at the following
function spaces:
\begin{equation}
\begin{split}
\textbf{V} &=
\overline{\Lambda_{c,\sigma}(\mathbb{H}^2(-a^2))}^{\|\cdot\|_{H^1(\mathbb{H}
^2(-a^2))}}, \\
\textbf{H} &=
\overline{\Lambda_{c,\sigma}(\mathbb{H}^2(-a^2))}^{\|\cdot\|_{L^2(\mathbb{H}
^2(-a^2))}},
\end{split}
\nonumber
\end{equation}
where $\Lambda_{c,\sigma}^1(\mathbb{H}^2(-a^2))$ is the space of all smooth,
compactly supported, divergence free $1$-forms on $\mathbb{H}^2(-a^2)$.
However, the work of the first and the second author \cite{CzubakChanRemarks}  indicates that the space
\begin{equation}
\widetilde{\textbf{V}} = \big \{  \phi \in H^1(\mathbb{H}^2(-a^2)) : \dd^* \phi
= 0  \big \}.
\nonumber
\end{equation}
  is much bigger than $\textbf{V}$.  Indeed, the
following orthogonal decomposition holds (this is the Helmholtz-Weyl decomposition, or can be thought of as a part of the Hodge decomposition, on the level of $H^1$ instead of the usual $L^2$)
\begin{equation}
\widetilde{\textbf{V}} = \textbf{V} \oplus \mathbb{F} ,
\nonumber
\end{equation}
where $\mathbb{F}$ stands for the space of all $L^2$-harmonic forms on
$\mathbb{H}^2(-a^2)$, that is,
\begin{equation}
\mathbb{F} = \big \{ \dd F \in L^2 (\mathbb{H}^2(-a^2)) : \Delta_H F = 0 \big
\}.
\end{equation}

From the above setting, we can now write down the most reasonable weak
formulation for the Cauchy problem of \eqref{fakeequation} as follows.

For any $v_0 \in \textbf{H}\oplus \mathbb{F}$, and external force $f \in L^2(0,T ;
\widetilde{\textbf{V}}') $, a weak solution to \eqref{fakeequation} which arises
from the initial datum $v_0$ is an element
\begin{equation}
v \in C^0([0,T] ; \textbf{H}\oplus \mathbb{F} ) \cap L^2(0,T ;
\widetilde{\textbf{V}}) ,
\nonumber
\end{equation}
which satisfies the following properties:
\begin{itemize}
 \item [I)] $\partial_t v \in L^2(0,T ; \widetilde{\textbf{V}}' ).$
 \item [II)] \begin{equation}\label{Weakhyperbolic}
\big < \partial_t v (t) , \phi  \big >_{\widetilde{\textbf{V}}'\otimes
\widetilde{\textbf{V}}} + \mu \int_{\mathbb{H}^2(-a^2)} g(\dd v(t) , \dd \phi )
\vol_{\mathbb{H}^2(-a^2)} =
\big < f(t) , \phi  \big >_{\widetilde{\textbf{V}}'\otimes
\widetilde{\textbf{V}}},
\end{equation}
holds for all $\phi \in \widetilde{\textbf{V}}$, and for almost every $t\in [0,T]$.
\item [III)]  $v(0) =v_0$.
\end{itemize}

Next, consider the following energy estimate, which is the natural analog of estimate \eqref{Euclideanapriori}
\begin{equation}\label{EnergyM}
\|v\|_{L^{\infty}(0,T ; L^2)}^2 + \mu \int_0^T
\int_{\mathbb{H}^2(-a^2)} \big | \nabla v (\tau )\big |^2
\vol_{\mathbb{H}^2(-a^2)} \dd \tau
\leq \|v_0\|_{L^2(\Omega)}^2 + \frac{1}{\mu} \|f\|_{L^2(0,T; \widetilde{\textbf{V}'})}^2.
\end{equation}
We now attempt to derive it from the weak formulation of the Stokes equation \eqref{fakeequation}.
Repeating the same type of Young's inequality estimate as we did in
\eqref{Younginequality}, now with \eqref{WeaklinearStokes} being replaced by
\eqref{Weakhyperbolic}, one gets instead
\begin{equation}\label{troubleestimate}
\begin{split}
\frac{1}{2} \frac{\dd }{\dd t} \|v(t)\|_{L^2(\mathbb{H}^2(-a^2))}^2 + \mu \|\dd v
(t)\|_{L^2(\mathbb{H}^2(-a^2))}^2 & =
\big < f(t) , v(t) \big >_{\widetilde{\textbf{V}}'\otimes \widetilde{\textbf{V}}
} \\
& \leq \frac{\epsilon}{2} \|\nabla v(t) \|_{L^2(\mathbb{H}^2(-a^2)}^2 +
\frac{1}{2 \epsilon} \|f(t)\|^2_{\widetilde{\textbf{V}}'} .
\end{split}
\end{equation}
The trouble in \eqref{troubleestimate} comes from the fact that one cannot use
$\mu \|\dd v(t)\|_{L^2(\mathbb{H}^2(-a^2))}^2$ to dominate the term
$\frac{\epsilon}{2} \|\nabla v(t)\|_{L^2(\mathbb{H}^2(-a^2))}^2$,
regardless of how small $\epsilon > 0 $ is. This is due to the non-equivalence
between $\|\dd \phi \|_{L^2(-a^2)}$ and $\|\nabla \phi
\|_{L^2(\mathbb{H}^2(-a^2))}$. If we choose a non-zero harmonic form $\phi = \dd
F \in L^2(\mathbb{H}^2(-a^2))$ with $F$ to be a harmonic function on
$\mathbb{H}^2(-a^2)$, then we will have $\|\nabla \dd F
\|_{L^2(\mathbb{H}^2(-a^2))} = a \|\dd F\|_{L^2(\mathbb{H}^2(-a^2))}> 0 $,
though we know $\|\dd \dd F\|_{L^2(\mathbb{H}^2(-a^2))} = 0$.

\begin{remark}
To fix this, one could be tempted to define a different norm on  $H^1(\mathbb{H}^2(-a^2))$ by $\|\dd v\|_{L^2(\mathbb{H}^2(-a^2))}$. However, then again considering a non-zero harmonic form $\phi=\dd
F \in L^2(\mathbb{H}^2(-a^2))$ with $F$ harmonic, we would get $\| \phi\|_{H^1(\mathbb{H}^2(-a^2))}=\|\dd \phi\|_{L^2(\mathbb{H}^2(-a^2))}=0$, but $\dd F\neq 0$, so this cannot be a norm on $H^1$.
\end{remark}

 While the above remarks do not constitute a proof that the estimate cannot be derived by any other means,
they suggest that
the non-equivalence
between $\|\dd \phi \|_{L^2(-a^2)}$ and $\|\nabla \phi
\|_{L^2(\mathbb{H}^2(-a^2))}$
creates immediate difficulties. In fact, we can prove that the estimate cannot be established, as the theorem below illustrates.

\begin{theorem}\label{secondtheorem}
There does not exist any absolute constant $C_0 >0$ for which the energy inequality
\begin{equation}\label{GEIingeneral}
\begin{split}
& \big \|v(t) \big \|_{L^{\infty}(0,T ;L^2(\mathbb{H}^2(-a^2)))}^2 + \int_0^T
\int_{\mathbb{H}^2(-a^2)} \big | \nabla v (\tau )\big |^2
\vol_{\mathbb{H}^2(-a^2)} \dd \tau \\
 &\quad \leq  C_0 \Big \{ \big \| v(0) \big \|_{L^2(\mathbb{H}^2(-a^2))}^2 + \big \| f  \big \|_{L^2 (0,T; \widetilde{\textbf{V}}' )}^2 \Big \}
 ,
\end{split}
\end{equation}
holds for any given terminal time $T> 0$, and any
$v\in C^0 ([0,T]; \textbf{H}\oplus \mathbb{F}) \cap L^2(0,T; \widetilde{\textbf{V}})$ satisfying the weak formulation \eqref{Weakhyperbolic}.
\end{theorem}

\begin{proof}
Assume towards contradiction that there was an absolute constant $C_0 > 0$, which depends only on the viscosity coefficient $\mu$ and the constant sectional curvature $-a^2$ of $\mathbb{H}^2(-a^2)$, such that \eqref{GEIingeneral} holds for any terminal time $T > 0$, and any $v\in C^0 ([0,T]; \textbf{H}\oplus \mathbb{F}) \cap L^2(0,T; \widetilde{V})$ satisfying the weak formulation \eqref{Weakhyperbolic}. To derive a contradiction, let
\begin{equation}\label{funnyv}
v(t) = \frac{4}{3} t^{\frac{3}{4}} \dd F ,
\end{equation}
where $\dd F \in L^2(\mathbb{H}^2(-a^2))$ with $F$ to a harmonic function on
$\mathbb{H}^2(-a^2)$. Then, we consider the following external forcing term
\begin{equation}\label{funnyF}
f(t) = t^{-\frac{1}{4}} \dd F.
\end{equation}
It is evident that $v \in C^0( [0,T] ; \textbf{H}\oplus \mathbb{F}) \cap
L^2(0,T ; \widetilde{\textbf{V}}  )$ and that
$f \in L^2(0,T ; \widetilde{\textbf{V}}' )$. It is equally evident that $v$
satisfies Properties I, II, and III, and hence is a weak solution to
\eqref{fakeequation} which arises from the initial datum $v(0) = 0$.

Now, the definition of $v$ at once gives
\begin{equation}\label{extraeasy1}
\big \| v \big \|_{L^{\infty} (0,T ; L^2(\mathbb{H}^2(-a^2)))}^2 = \frac{16}{9} T^{\frac{3}{2}} \big \| \dd F \big \|_{L^2(\mathbb{H}^2(-a^2))}^2 .
\end{equation}
It follows also from a direct computation that
\begin{equation}\label{extraeasy2}
\int_0^T \big \| \nabla v(t) \big \|_{L^2(\mathbb{H}^2(-a^2))}^2 \dd t = \frac{32}{45} a^2 T^{\frac{5}{2}} \big \| \dd F \big \|_{L^2(\mathbb{H}^2(-a^2))}^2 .
\end{equation}
It is equally obvious that the following relation holds for any $t \in [0,T]$
\begin{equation*}
\big \| f(t) \big \|_{\widetilde{\textbf{V}}' }^2 \leq t^{-\frac{1}{2}} \big \| \dd F \big \|_{L^2(\mathbb{H}^2(-a^2))}^2 ,
\end{equation*}
from which we deduce
\begin{equation}\label{extraeasy3}
\int_0^T \big \| f(t) \big \|_{\widetilde{\textbf{V}}' }^2 \dd t \leq 2 T^{\frac{1}{2}} \big \| \dd F \big \|_{L^2(\mathbb{H}^2(-a^2))}^2 .
\end{equation}
Since we insist that global energy estimate \eqref{GEIingeneral} is valid, by combining \eqref{extraeasy1}, \eqref{extraeasy2}, and \eqref{extraeasy3}, it would follow from \eqref{GEIingeneral} that the following estimate would hold for any $T > 0$ ,
\begin{equation}
\frac{16}{9} T^{\frac{3}{2}} + \frac{32}{45} T^{\frac{5}{2}} a^2 \leq 2 C_0 T^{\frac{1}{2}} ,
\end{equation}
which however is clearly absurd. So, a contradiction has been derived.
\end{proof}

Next, let us mention very briefly what could happen if we consider the following version of the incompressible Navier-Stokes equation instead.

\begin{equation}\label{fakeNSequation}
\left\{
\begin{split}
\partial_t v -\mu\Delta_H v + \nabla_v v + d P & = f , \\
\dd^* v & = 0.
\end{split}
\right.
\end{equation}

Observe that by taking $v$ to be the same one as in \eqref{funnyv} and taking the following function.
\begin{equation}\label{funnypressure}
P = -\frac{1}{2} \Big ( \frac{4}{3} t^{\frac{3}{4}} \Big )^2 \big | \dd F \big |^2 ,
\end{equation}
it follows that $v$ is a smooth solution to \eqref{fakeNSequation} on $[0,\infty )\times \mathbb{H}^2(-a^2)$, with $f$ and $P$ to be specified in \eqref{funnyF} and \eqref{funnypressure} respectively. This means that exactly the same argument as given in the proof of Theorem \ref{secondtheorem} will give us the following result.

\begin{theorem}
There does not exist any absolute constant $C_0 > 0$ such that the a priori energy estimate \eqref{GEIingeneral} holds for any terminal time $T > 0$ and any classical solution $v$  to \eqref{fakeNSequation} on $[0,T]\times \mathbb{H}^2(-a^2)$ which satisfies $v \in L^{\infty} (0,T ; L^2(\mathbb{H}^2(-a^2))) \cap L^2(0,T ; H^1(\mathbb{H}^2(-a^2)))$ .
\end{theorem}

\begin{remark}
In the setting of the round sphere $S^2(\frac{1}{a})$ with radius $\frac{1}{a}$,
the derivation of energy-type apriori estimate for a solution directly from
\eqref{fakeequation} works well however. This is due to the equivalence between
$\|\dd \phi\|_{L^2(S^2(\frac{1}{a}))}$ and $\|\nabla \phi
\|_{L^2(S^2(\frac{1}{a}))}$ for all $\phi \in H^1(S^2(\frac{1}{a}))$ with $\dd^*
\phi = 0 $. In this setting, the Bochner-Weitzenbock's formula reads
$\nabla^*\nabla \phi = (-\Delta_H) \phi - a^2 \phi$, from which it follows,
through integration by parts, that the following relation holds for all $\phi
\in H^1(S^2(\frac{1}{a}))$ satisfying $\dd^* \phi = 0$.
\begin{equation}
\begin{split}
\|\dd \phi \|_{L^2(S^2(\frac{1}{a}))}^2 & = \|\nabla \phi
\|_{L^2(S^2(\frac{1}{a}))}^2 + a^2 \| \phi \|_{L^2(S^2(\frac{1}{a}))}^2 \\
& \geq  \|\nabla \phi \|_{L^2(S^2(\frac{1}{a}))}^2 .
\end{split}
\nonumber
\end{equation}
On the other hand, we have the following expression for $\dd$ in  terms
of $\nabla$
\begin{equation}\label{goodexpressions}
\dd = \eta^{\alpha}\wedge \nabla_{e_{\alpha}} ,
\end{equation}
with $\{ e_1, e_2 \}$ and $\{ \eta^1 , \eta^2 \}$ positively oriented
orthonormal frame and dual frame on $S^2(\frac{1}{a})$.  Equivalently, we can write in coordinates
\begin{equation}\label{goodexpressions2}
\begin{split}
\dd v & = \frac 12(\partial_i v_j -\partial_j v_i)\dd x^i \wedge \dd x^j=  \frac 12(\nabla_i v_j -\nabla_j v_i)\dd x^i \wedge \dd x^j\quad\mbox{since} \ \Gamma^l_{ij}=\Gamma^l_{ji},
\end{split}
\end{equation}
where we sum over repeated indices.
Then \eqref{goodexpressions}  or  \eqref{goodexpressions2} gives
\begin{equation}
\begin{split}
\big \| \dd v \big \|_{L^2(\mathbb{H}^2(-a^2))} & \leq 2^{\frac{1}{2}} \big \| \nabla v  \big \|_{L^2(\mathbb{H}^2(-a^2))}.
\end{split}
\end{equation}

This successful attempt in the
spherical case should not be viewed as evidence in favor of the Hodge Laplacian
as an admissible candidate for the viscosity operator. This is because the
naturality and correctness of a right choice of the operator representing
viscosity should be equally valid for (and hence consistent with the settings
of) all possible Riemannian manifolds.
\end{remark}

Next, we consider \eqref{fakeequation} on $\mathbb{H}^2(-a^2)$  from the
perspective of a
 global energy equality. According to standard theory of linear parabolic
equations, one expects that the following global energy equality should be
derived as a consequence of the governing parabolic equation
\eqref{fakeequation}, and its weak formulation.
\begin{equation}\label{GEE}
\begin{split}
& \frac{1}{2}\|v(t)\|_{L^2(\mathbb{H}^2(-a^2))}^2 + \mu \int_0^t
\int_{\mathbb{H}^2(-a^2)} \big | \nabla v (\tau )\big |^2
\vol_{\mathbb{H}^2(-a^2)} \dd \tau \\
 &\quad= \frac{1}{2}\|v(0)\|_{L^2(\mathbb{H}^2(-a^2))}^2 + \int_0^t \big < f(t) ,
v(\tau) \big >_{\widetilde{\textbf{V}}'\otimes \widetilde{\textbf{V}} } \dd \tau
.
\end{split}
\end{equation}
\begin{remark}
   In the setting of $\mathbb{R}^2$ or $\mathbb{R}^3$, the global energy
equality in exactly the same form as that of \eqref{GEE} (simply with
$\mathbb{H}^2(-a^2)$ replaced by $\mathbb{R}^n$ or its bounded sub-domains with
smooth boundaries) is indeed derived directly from \eqref{Weakhyperbolic} (see,
for instance, chapter 4 in \cite{SereginBook}). This basic fact is exactly what
motivates us to ask for the validity of \eqref{GEE} in the setting of
$\mathbb{H}^2(-a^2)$.
\end{remark}
So, can one derive \eqref{GEE} directly from the weak formulation
\eqref{Weakhyperbolic}? The following theorem provides a counterexample.
\begin{theorem}\label{FirstTheorem}
The energy equality \eqref{GEE} cannot be directly derived from the weak formulation
\eqref{Weakhyperbolic}.
\end{theorem}
\begin{proof}
We again look at the velocity $v$ which is specified in \eqref{funnyv}, as well as the external forcing term $f$ which is specified in \eqref{funnyF}. Again, such a pair $(v, f)$ satisfies the weak formulation \eqref{Weakhyperbolic}.
 
However, a direct computation gives
\begin{equation}\label{assistance}
\frac{1}{2} \|v(t)\|_{L^2(\mathbb{H}^2(-a^2))}^2 = \frac{8}{9} t^{\frac{3}{2}}
\|\dd F\|_{L^2(\mathbb{H}^2(-a^2))}^2
= \int_0^t \big < f(\tau ) , v (t) \big >_{\widetilde{\textbf{V}}' \otimes
\widetilde{\textbf{V}}} .
\nonumber
\end{equation}
Again $v(0)= 0$. So, if we insist that the global energy equality
\eqref{GEE} would hold for $v$ as given in \eqref{funnyv}, this would lead to
the following
\begin{equation}
\int_0^T \int_{\mathbb{H}^2(-a^2)} \big | \nabla v(t) \big |^2
\vol_{\mathbb{H}^2(-a^2)} \dd \tau = 0 ,
\nonumber
\end{equation}
which in turn would lead to the contradictory statement that $v$ is identically zero
on $[0,T]\times \mathbb{H}^2(-a^2)$.
\end{proof}

Thus, one cannot choose the Hodge Laplacian
as the viscosity operator in the general setting of a Riemannian manifold if
one hopes to obtain a global energy inequality or equality. Given that a global energy inequality
is one of the cornerstones of the existence theory for the Navier-Stokes equations, and that the global energy equality is a natural consequence,
it seems very
likely that a successful theory of the Navier-Stokes equations on Riemannian
manifolds cannot
rely on such a choice of the viscosity operator.

\begin{remark}
We remark that an existence theory for the Navier-Stokes equations on manifolds using the Hodge-Laplacian has been
developed by Ebin and Marsden in \cite{EbinMarsden}. Our counter-example, however, is not at odds
with \cite{EbinMarsden} because there the authors treat the case of compact manifolds, whereas here
we construct a counterexample on the (non-compact) hyperbolic plane.  In addition, the theory developed in \cite{EbinMarsden} is a short-time result, which requires the initial data to be in $H^s$ for $s>\frac n2+5$, whereas the discussion here pertains to the global in time Leray-Hopf solutions for which one only assumes initial data in $L^2$.
\end{remark}

\section{Restriction arguments\label{section_restriction}}

Faced with more than one choice for the viscosity operator on a Riemannian
manifold,
one may attempt to settle the question upon examination of the Euclidean case,
as follows.
Since all natural choices for the viscosity operator agree in $\R^3$, we can
try to obtain
its correct form on manifolds by means of analyzing the restriction of $\Delta
v$ to the two-sphere $S^2$,
where $v$ is a (divergence-free) vector field in $\R^3$. While one could do
this for any embedded submanifold and
any dimension, the case of $S^2$ is natural in light of the well-known formula
\begin{gather}
\Delta = \partial^2_r + \frac{2}{r} \partial_r + \frac{1}{r^2} \Delta_{S^2}
\label{decomposition_sphere_functions}
\end{gather}
for the Laplacian in spherical coordinates $(r,\varphi,\theta)$ acting on
functions.  We show that if we start with the divergence free vector field on $\R^3$, taking the Hodge Laplacian on $\R^3$, and then restricting it to a sphere $S^2$, we indeed 
can obtain
\beq\label{restr}
(-\Delta_H v)\large|_{S^2}=-2v\large|_{S^2}-\Delta_{S^2}(v\large|_{S^2})=-2\Ric_{S^2}(v\large|_{S^2})-\Delta_{S^2}(v\large|_{S^2}),
\ee 
where $-\Delta_{S^2}$ denotes now the Hodge Laplacian on $S^2$.  Then \eqref{divdef} gives another indication that the deformation tensor is the correct operator to use.

To derive \eqref{restr},  consider polar coordinates
$(r,\varphi,\theta)$  and let
\begin{gather}
e_1 = \partial_r,\,
e_2 = \frac{1}{r} \partial_\varphi, \,
e_3 = \frac{1}{r\sin\varphi} \partial_\theta,
\nonumber
\end{gather}
so that a dual basis is given by
\begin{gather}
(e_1)^* \equiv  e^1 = dr, \,
(e_2)^* \equiv e^2 = r \, d\varphi, \,
(e_3)^* \equiv e^3 = r \sin\varphi\, d\theta.
\nonumber
\end{gather}
One easily computes
\begin{gather}
d\theta \wedge d\varphi = \frac{1}{r^2 \sin\varphi} e^3 \wedge e^2,
\nonumber
\end{gather}
\begin{gather}
d\varphi \wedge dr = \frac{1}{r} e^2 \wedge e^1,
\nonumber
\end{gather}
\begin{gather}
d\theta \wedge dr = \frac{1}{r \sin\varphi} e^3 \wedge e^1,
\nonumber
\end{gather}
and, of course, $\vol = e^1 \wedge e^2 \wedge e^3$.
Then recalling the definition of the Hodge $\ast$ operator
\[
\alpha \wedge \ast \beta =g(\alpha, \beta) \vol,
\]
we compute
%
\beq\label{ast2}
\ast(dr\wedge d\varphi)=\sin \varphi d\theta,\quad \ast(dr\wedge d\theta)=-\frac{1}{\sin\varphi}d\varphi,\quad \ast (d\varphi\wedge d\theta)=\frac{1}{r^2\sin\varphi}dr.
\ee

Now let $v_\sharp$ be a vector field on $\R^3$ that is independent of $r$
\[
v_\sharp=v_\varphi \partial_\varphi +v_\theta \partial_\theta,
\]
so it is also a vector field on $S^2.$  Assume that $v$ is divergence free on $\R^3$
\beq\label{divfree}
\partial_\varphi(\sin \varphi v_\varphi)+\sin \varphi \partial_\theta v_\theta=0.
\ee
Observe this is also equivalent to $v_\sharp$ being divergence free on $S^2$.
Lowering the index, we get a $1$-form 
\begin{gather}
v = r^2v_\varphi \, d\varphi + r^2\sin^2\varphi v_\theta \, d\theta.
\nonumber
\end{gather}
We now can use the previous identities and definitions to find
\begin{align}
\begin{split}
* dv & = \frac{\partial_r (r^2v_\varphi)}{r} e^3
- \frac{ \partial_r (r^2\sin^2\varphi v_\theta)}{r \sin\varphi} e^2
- \frac{\partial_\theta (r^2v_\varphi) - \partial_\varphi (r^2\sin^2\varphi v_\theta)}{r^2 \sin\varphi}
e^1
\\
&
=
\partial_r (r^2v_\varphi)  \sin\varphi \, d\theta
-\frac{\partial_r (r^2\sin^2\varphi v_\theta)}{ \sin\varphi} \, d\varphi
+\frac{\partial_\varphi (r^2\sin^2\varphi v_\theta)-\partial_\theta (r^2v_\varphi) }{r^2 \sin\varphi}
\, dr\\
&
=
2r v_\varphi  \sin\varphi \, d\theta
-2r \sin \varphi v_\theta \, d\varphi
+\frac{\partial_\varphi (\sin^2\varphi v_\theta)-\partial_\theta v_\varphi }{ \sin\varphi}
\, dr.
\end{split}
\nonumber
\end{align}

Next recall that
\begin{gather}
d^* = (-1)^{n(p+1) + 1} * d *,
\nonumber
\end{gather}
where $p$ is the degree of the form $d^\ast$ is acting on,
so on two forms in  $\R^3$,
\begin{gather}
d^* = *d* .
\nonumber
\end{gather}
Since $v_\sharp$ is divergence-free, we have $d^* v = 0$, so that 
\[
-\Delta_H v  = d d^*v + d^* dv = d^* dv = *d*dv.
\]
Then computing further  

\begin{align*}
d\ast dv&=
\left(\partial_r(2r v_\varphi  \sin\varphi )-\partial_\theta (\frac{\partial_\varphi (\sin^2\varphi v_\theta)-\partial_\theta v_\varphi }{ \sin\varphi}) \right )dr\wedge d\theta\\
&\quad -\left( \partial_r(2r \sin \varphi v_\theta)+\partial_\varphi(  \frac{\partial_\varphi (\sin^2\varphi v_\theta)-\partial_\theta v_\varphi }{ \sin\varphi}
  ) \right)dr\wedge d\varphi\\
 &\quad
+\left(\partial_\varphi(2r v_\varphi  \sin\varphi )+\partial_\theta(2r \sin \varphi v_\theta) \right) d\varphi \wedge d \theta\\
&=\left(2v_\varphi  \sin\varphi -\frac{\partial_\theta (\partial_\varphi (\sin^2\varphi v_\theta)-\partial_\theta v_\varphi )}{ \sin\varphi}\right) dr\wedge d\theta\\
&\quad-\left(2 \sin \varphi v_\theta+\partial_\varphi(  \frac{\partial_\varphi (\sin^2\varphi v_\theta)-\partial_\theta v_\varphi }{ \sin\varphi}
  )\right) dr\wedge d\varphi,
\end{align*}
by the divergence free condition \eqref{divfree}.
Taking $\ast$ by using \eqref{ast2}, we get
\begin{align*}
-\Delta_H v  = *d*dv&
 =-2v_\varphi d\varphi +\frac{\partial_\theta (\partial_\varphi (\sin^2\varphi v_\theta)-\partial_\theta v_\varphi )}{ \sin^2\varphi} d\varphi\\
&\quad-2\sin^2 \varphi v_\theta d\theta-\partial_\varphi(  \frac{\partial_\varphi (\sin^2\varphi v_\theta)-\partial_\theta v_\varphi }{ \sin\varphi}
  ) (\sin \varphi d\theta)\\
 &  =
 -2v\large|_{S^2}+\frac{\partial_\theta (\partial_\varphi (\sin^2\varphi v_\theta)-\partial_\theta v_\varphi )}{ \sin^2\varphi} d\varphi
 -\partial_\varphi(  \frac{\partial_\varphi (\sin^2\varphi v_\theta)-\partial_\theta v_\varphi }{ \sin\varphi}
  ) (\sin \varphi d\theta) \\
&  = -2v\large|_{S^2}- 
  \Delta_{S^2} v \large|_{S^2}
\end{align*}
as needed.

\section{Non-relativistic limit\label{section_nr_limit}}

Another way we can investigate the formulation of the Navier-Stokes equations
on Riemannian manifolds is via the non-relativistic limit of the relativistic Navier-Stokes
equations, as follows.
Starting with a relativistic form of the momentum
equation, we can
take the non-relativistic limit characterized by fluid velocities very small
compared to the speed of light.
But differently than the usual
non-relativistic limit in the general theory of relativity, where it is assumed
that the metric converges to
the Minkowski metric, with the metric induced on $\{ t = \text{ constant}\}$
hypersurfaces being,
therefore, the Euclidean one, we can consider the situation where the metric
on
$\{ t = \text{ constant}\}$ hypersurfaces converges to an arbitrary Riemannian
metric. This procedure leads to an evolution equation for the (non-relativistic) velocity
and pressure. Considering further the incompressible limit, we arrive at a candidate
for the Navier-Stokes equations on Riemannian manifolds. Because  general relativity is a
more accurate description of nature than classical mechanics, equations that are obtained
from the former are well-motivated from a physical perspective.

An immediate difficulty to this approach is that 
the correct
formulation of relativistic viscous fluids is not known. The most
direct attempts of generalizing  (\ref{NSR}) to relativity,
introduced by Eckart  \cite{EckartViscous} and Landau and Lifshitz
 \cite{LandauLifshitzFluids},
 lead to many
patahologies,
 including a breakdown of causality \cite{Hiscock_Lindblom_instability_1985,
Hiscock_Lindblom_pathologies_1988,PichonViscous}. Therefore, such theories
cannot be taken as a physically sensible model for the study of the non-relativistic limit.

Alternatives to
 \cite{EckartViscous,LandauLifshitzFluids} have been proposed
by Lichnerowicz \cite{LichnerowiczBookGR}, Choquet-Bruhat
\cite{ChoquetBruhatGRBook},
and Freist\"uhler and Temple \cite{TempleViscous}. Each of these has been shown
to yield a satisfactory theory of relativistic viscous fluids under different
assumptions
\cite{ChoquetBruhatGRBook, DisconziCzubakNonzero,
DisconziViscousFluidsNonlinearity,
TempleViscous}. Although these results fall short of covering all situations of
physical interest,
and the matter of how to correctly formulate relativistic viscous phenomena
remains open, they provide a basis for the study of the non-relativistic limit, and 
a motivation for a particular form of the Navier-Stokes equations on Riemannian manifolds.

What is important from our perspective is that the equations introduced 
in 
\cite{ChoquetBruhatGRBook, TempleViscous, LichnerowiczBookGR},
all lead to the same form of the 
(non-relativistic) Navier-Stokes on Riemannian manifolds, namely, (\ref{NS}).

In a certain sense, this is not a surprise: the relativistic equations are obtained
as 
\begin{gather}
\nabla^\al T_{\al\be} = 0
\label{div_T_rel}
\end{gather}
 ($\al,\be=0,\dots, 3$), where $T_{\al\be}$ is a symmetric
two-tensor. Thus, if $T_{\al\be}$ depends on first-derivatives of the velocity, it will
in general contain the term $\nabla_\al u_\be +\nabla_\be u_\al$, where $u$ is the 
fluid's four-velocity (which is the appropriate relativistic notion of velocity). This is indeed
the case in all the relativistic theories mentioned above. It can be showed
(see, e.g., \cite[Chapter 6]{RezzollaZanottiBookRelHydro} or \cite{DisconziKephartScherrerNew})
 that 
after taking the non-relativistic limit of (\ref{div_T_rel}), the term 
$\nabla_\al u_\be +\nabla_\be u_\al$ will produce $\nabla^i(\nabla_i v_j + \nabla_j v_i)$, where
$v$ is the (classical) fluid's velocity ($i,j=1,2,3$). We see that this gives precisely
the viscosity operator in (\ref{NS}) in view of (\ref{divdef}).

\begin{remark}
There exist other approaches to viscosity in relativity that are based on 
ideas of relativistic extended irreversible
thermodynamics
\cite{JouetallBook, MuellerRuggeriBook}, most notably
the Mueller-Israel-Stewart theory \cite{MIS-2, MIS-3, MIS-5, MIS-6, MIS-1,
MIS-4}
and the divergence-type theories \cite{GerochLindblomDivergenceType, GerochLindblomCausal,LiuMullerRuggeri-RelThermoGases}. However,
in these theories 
the viscous contributions to $T_{\al\be}$ are not given in terms
of $u$ and
the other thermodynamic variables, but rather are treated as additional variable
in the problem. Therefore, they are not suitable for determining the form 
of the Navier-Stokes equations on manifolds via a non-relativistic limit.
\end{remark}

\section*{Acknowledgments}
The authors would like to thank the referee for the careful reading of the manuscript and the insightful comments and suggestions.  We also would like to thank Steven Preston and Doug Wright for independently suggesting the restriction argument to the sphere, and we thank the referee for recommending to consider the vector field instead of the $1-$form.

Chi Hin Chan is partially supported by a
grant from the National Science Council of Taiwan (NSC 101-2115-M-009-016-MY2).
Magdalena Czubak is partially supported by a grant from the Simons Foundation \# 246255.
Marcelo M. Disconzi is partially supported by NSF grant 1305705.
This work has been partially supported by Vanderbilt International Research Grant.

\bibliographystyle{plain}
\bibliography{References_Marcelo.bib}

\end{document}